\documentclass{amsart}

\pdfoutput = 1

\usepackage{amsfonts}
\usepackage{amssymb}
\usepackage{amsthm}
\usepackage{amsmath}
\usepackage[all]{xy}
\usepackage{hyperref}

\DeclareMathOperator{\Hom}{Hom}
\DeclareMathOperator{\Ext}{Ext}

\DeclareMathOperator{\gr}{gr}
\DeclareMathOperator{\img}{im}
\DeclareMathOperator{\Cob}{Cob}
\DeclareMathOperator{\spn}{span}
\DeclareMathOperator{\E}{E}
\newcommand{\ktwo}{\mathcal{K}_2}
\newcommand{\B}{\mathcal{B}}
\newcommand{\T}{\mathbb{T}}
\newcommand{\M}{\mathcal{M}}
\newcommand{\R}{\mathcal{R}}
\newcommand{\G}{\mathcal{G}}
\newcommand{\bbk}{\mathbb{K}}
\newcommand{\Gr}{\mathrm{Gr}}
\newcommand{\EGr}{\E_\Gr}

\newtheorem{theorem}{Theorem}[section]

\newtheorem{lemma}[theorem]{Lemma}
\newtheorem{proposition}[theorem]{Proposition}

\theoremstyle{definition} 
\newtheorem{definition}[theorem]{Definition}
\newtheorem{example}[theorem]{Example}

\begin{document}

\title{Generalized Koszul properties for augmented algebras}

\author[Christopher Phan]{}



\maketitle

\begin{center}
\vskip -4 mm
{\sc Christopher Phan}
\vskip 5 mm
Department of Mathematics\\
1222 University of Oregon\\
Eugene, Oregon 97403\\
{\tt cphan@uoregon.edu}
\end{center}

\vskip 5 mm

\begin{abstract}
Under certain conditions, a filtration on an augmented algebra $A$ admits a related filtration on the Yoneda algebra $\E(A) := \Ext_A(\bbk, \bbk)$. We show that there exists a bigraded algebra monomorphism $\gr \E(A) \hookrightarrow \EGr(\gr A)$, where $\EGr(\gr A)$ is the graded Yoneda algebra of $\gr A$. This monomorphism can be applied in the case where $A$ is connected graded to determine that $A$ has the $\ktwo$ property recently introduced by Cassidy and Shelton.
\end{abstract}


\section{Introduction}
In this paper, we use filtrations to study certain homological properties of augmented algebras. We generalize a similar recently-studied homological property of graded algebras. Throughout, if a $\bbk$-algebra $A$ (where $\bbk$ is a field) is graded by a monoid $\M$ with identity element $e$, we denote by $\Ext_{\Gr}$ the derived functor of the $\M$-graded $\Hom$ functor 
\begin{equation*}
\Hom_{\Gr}(M,N) := \bigoplus_{\alpha \in \M} \Hom_{\Gr}(M, N)_\alpha,
\end{equation*}
where $\Hom_{\Gr}(M, N)_\alpha = \hom_A(M(\alpha),N)$,
$M(\alpha)_\beta := M_{\alpha\beta}$, and $\hom_A(M,N)$ contains only degree-$e$ homomorphisms. A connected-graded algebra $A$ is called {\bf Koszul} if its (graded) Yoneda algebra $\EGr(A) := \Ext_{\Gr}(\bbk, \bbk)$ is generated as a $\bbk$-algebra by $\EGr^1(A)$. (Throughout, we assume connected-graded algebras are finitely generated and finitely related.) Our goal is to study a generalization of Koszul introduced by Cassidy and Shelton \cite{cassidysheltonK2alg06}:
\begin{definition}\label{d:ktwograded} A connected-graded algebra $A$ is $\mathcal{K}_2$ if $\EGr^1(A)$ and $\EGr^2(A)$ generate $\EGr(A)$ as a $\bbk$-algebra.
\end{definition}

The Koszul algebras are exactly the quadratic $\ktwo$ algebras. More generally, an algebra $A$ is in the class of $N$-Koszul algebras introduced by Berger \cite{berger01} if and only if $A$ is $\ktwo$ and has only degree-$N$ relations \cite[Theorem 4.1]{green04}\cite[Corollary 4.6]{cassidysheltonK2alg06}. (Unfortunately, the term {\it $N$-Koszul} has obtained two incompatible meanings. The meaning used here is different than that found in \cite{polishquadalg}.) However, $\ktwo$ algebras can have relations of several different degrees.

Our goal is generalize further to augmented algebras and to relate the graded case to the augmented case via filtrations. Throughout, we suppose that $A$ is an augmented algebra over a field $\bbk$, i.e., $A = A_+ \oplus \bbk \cdot 1$ for $A_+ \triangleleft A$. (The augmentation is then $\varepsilon: A \twoheadrightarrow \bbk$).  Furthermore, we suppose that a totally-ordered monoid $\M$ (with identity element $e$) filters $A$ so that 
\begin{enumerate}
\item $\bigcup_\alpha F_\alpha A = A$;
\item $F_\alpha A = \bbk \oplus F_\alpha A_+$, where $F_\alpha A_+ := F_\alpha A \cap A_+$;
\item $F_e A = \bbk$ and $F_\alpha A_+ \not = 0$ when $\alpha > e$; and
\item\label{i:filtrationdimcondo} $\dim F_\alpha / F_{s(\alpha, -1)}A < \infty$ for all $\alpha > e$.
\end{enumerate}

We use $\E(A)$ to denote the Yoneda algebra $\Ext_A(\bbk, \bbk)$, the cohomology of the cobar complex $\Cob(A) := \Hom_\bbk(A_+^{\otimes \bullet}, \bbk)$, where $\Hom$ is the functor yielding \emph{all} $A$-module homomorphisms. The complex $\Cob(A)$ has an $\M$-filtration $F_\alpha \Cob(A)$ (see Definition \ref{d:cobar}) which induces a filtration $F_\alpha \E^n(A)$ and associated graded algebra $\gr^F \E(A)$. Also, the filtration on $A$ yields the associated graded algebra $\gr^F A$ (graded by $\M$); we set $(\gr^F A)_+ := \bigoplus_{\alpha > e} (\gr^F A)_\alpha$. The algebra $\gr^F A$ is augmented by $\gr^F A = \bbk \oplus (\gr^F A)_+$. 

\begin{theorem}\label{t:algmonomorphism}
There is a bigraded (with respect to the cohomological and $\M$ gradings) algebra monomorphism
\begin{equation*}
\Lambda: \gr^F \E(A) \hookrightarrow \EGr(\gr^F A).
\end{equation*}
\end{theorem}

We make the following generalization of $\ktwo$ to this broader category of algebras:
\begin{definition}\label{d:ktwoungraded}
An augmented algebra $A$ is $\ktwo$ if $\E^1(A)$ and $\E^2(A)$ generate $\E(A)$ as a $\bbk$-algebra.
\end{definition} 
We can then connect the theory of connected-graded (finitely-related) algebras and ungraded algebras with the following, to be proved in Section \ref{s:monofilt}:
\begin{lemma}\label{t:findimequal}
For a connected-graded algebra $A$, $\E^m(A) = \EGr^m(A)$ if an only if $\dim \EGr^m(A) < \infty$. Consequently, if a connected-graded algebra $A$ is $\ktwo$ in the sense of Definition \ref{d:ktwoungraded}, then $A$ is $\ktwo$ in the sense of Definition \ref{d:ktwograded}.
\end{lemma}

Our primary goal was to develop a technique for transferring the $\ktwo$ property from $\gr^F A$ to $A$. As we will see, it is often much easier to prove that $\gr^F A$ is $\ktwo$.
\begin{theorem}\label{t:kmpush}
If $\EGr^1(\gr^F A)$ and $\EGr^2(\gr^F A)$ are finite dimensional and generate $\EGr(\gr^F A)$, and $\Lambda^1$ and $\Lambda^2$ are surjective, then $A$ is $\ktwo$.
\end{theorem}

This theorem captures a more specific situation involving connected-graded algebras. Every connected-graded algebra $A$ with $n$ generators is a factor of the free algebra $\bbk\left<x_1, \dots, x_n\right>$. Thus, the monomials in $x_1, \dots, x_n$ are a totally-ordered (noncommutative) monomial (under the degree-lexicographical order), and so provide a filtration $F$ on $A$. The following is well-known (see, for example, \cite[Theorem IV.3.1]{polishquadalg}):
\begin{theorem}\label{t:PBWalgs} If $A$ is a connected-graded quadratic algebra and $\gr^F A$ is also quadratic, then $A$ is Koszul.
\end{theorem} An algebra $A$ which meets the hypotheses of Theorem \ref{t:PBWalgs} is called a {\bf Poincar\'e--Birkhoff--Witt} algebra.
Setting $I := \ker(\bbk\left<x_1, \dots, x_n\right> \twoheadrightarrow A)$, we say a Gr\"obner basis $\mathcal{G}$ for $I$ is essential if its elements generate $I$ in a certain minimal manner (see Definition \ref{d:essgrob}). The following $\ktwo$ analogue of Theorem \ref{t:PBWalgs} was the original goal of this research.
\begin{theorem} \label{t:genPBWa}
If $I$ has an essential Gr\"obner basis and $\gr^F A$ is $\ktwo$, then $A$ is $\ktwo$ as well.
\end{theorem}

The algebra $\gr^F A$ will be a monomial connected-graded algebra. Cassidy and Shelton have provided an algorithm that determines whether a monomial connected-graded algebra is $\ktwo$ \cite[Theorem 5.3]{cassidysheltonK2alg06}. 

In Section \ref{s:algmonomorphism}, we prove Theorems \ref{t:algmonomorphism} and \ref{t:kmpush}, which involves relating the cobar complexes $\Cob(A)$ and $\Cob(\gr^F A)$, and constructing the map $\Lambda$. In Section \ref{s:monofilt}, we consider the case where $A$ is a connected-graded algebra, and connect the existence of an essential Gr\"obner basis to the surjectivity of $\Lambda_2$, proving Theorem \ref{t:genPBWa}. (The surjectivity of $\Lambda_1$ is automatic in the connected-graded case.) We connect surjectivity of $\Lambda^2$ with the existence of a special Gr\"obner basis for $\ker (\bbk\left<x_1, \dots, x_n\right> \twoheadrightarrow A)$. In Section \ref{s:anticommface}, we use the results from Section \ref{s:monofilt} to prove that some anticommutative analogues of face rings are $\ktwo$. 


\section{Bigraded algebra monomorphism $\Lambda: \gr^F \E(A) \hookrightarrow \EGr(\gr^F A)$} 
\label{s:algmonomorphism}

In this section, $A$ denotes an augmented algebra filtered by a totally-ordered monoid $\M$ as specified above. (Note that $\M$ need not be commutative.) We prove Theorems \ref{t:algmonomorphism} and \ref{t:kmpush}.
We shall use the following notation:
\begin{definition}
For $\alpha \in \M$, denote the monomial appearing $r$ steps after $\alpha$ in the total order by $s(\alpha, r)$.
\end{definition} 
We begin by setting detailed notation for the cobar complex and its associated filtration.

\begin{definition} \label{d:cobar} Let $d: A_+^{\otimes n} \rightarrow A_+^{\otimes n-1}$ via
\begin{equation*}
d(a_1 \otimes \cdots \otimes a_n) := 
\sum_{i = 1}^{n-1} (-1)^i a_1 \otimes \cdots \otimes a_ia_{i+1} \otimes \cdots \otimes a_n.
\end{equation*}
This makes $A_+^{\otimes \bullet}$ into a chain complex. We filter this complex by setting
\begin{equation*}
F_\alpha A_+^{\otimes n} := \sum_{\substack{\alpha_1 \cdots \alpha_n < \alpha \\
\alpha_i > e}} F_{\alpha_1} A_+ \otimes \cdots \otimes F_{\alpha_n} A_+.
\end{equation*}

The {\bf cobar complex} is the co-chain complex dual to $A_+^\bullet$, defined via
\begin{equation*}
\Cob^n(A) := \Hom_\bbk \left(A_+^{\otimes n}, \bbk\right)
\end{equation*}
with the dual differential, which we denote $\partial$.

We put a decreasing filtration on $\Cob(A)$ by setting
\begin{equation*}
F_\alpha \Cob^n(A) := \left\{f: A_+^{\otimes n} \rightarrow \bbk \, | \, F_{s(\alpha,-1)} A_+^{\otimes n} \subset \ker f\right\}.
\end{equation*}

If $B$ is an algebra graded by a totally-ordered monoid $\M$ with identity element $e$, we similarly define $\Cob_{\Gr}(B) := \Hom_{\Gr}\left(B_+^{\otimes \bullet}, \bbk\right)$, where $B_+ = \sum_{\alpha > e} B_\alpha$.

The cup product multiplication in a cobar complex, graded cobar complex, or Yoneda algebra will be denoted by $\smallsmile$.
\end{definition}
Throughout, we will denote $\Hom_\bbk(V, \bbk) =: V^\vee$. We first relate $\Cob^\bullet_{\Gr}(\gr^F A)$ to the cobar complex of $A$.

\begin{proposition} \label{t:dgaisom}
There is a differential-graded algebra isomorphism 
\begin{equation*}
\gr^F \Cob^\bullet(A) \simeq \Cob_{\Gr}^\bullet(\gr^F A).
\end{equation*}
\end{proposition}

The proof of Proposition \ref{t:dgaisom} will follow after two lemmas. 
Let us fix a $\bbk$-basis $\R = \coprod_{\alpha \in \M} \R_\alpha$ for $A$ such that: 
\begin{enumerate}
\item $\bigcup_{\beta \leq \alpha} \R_\beta$ is a basis for $F_\alpha A$.
\item $\R_\alpha \subset F_\alpha A_+$ for $\alpha > e$.
\end{enumerate}
Then $\{r + F_{s(\alpha, -1)} A_+ : r \in \R_\alpha\}$ is a basis for $F_\alpha A_+ / F_{s(\alpha, -1)} A_+$.

For readability, we set $(((\gr^F A)_+)^{\otimes n})_\alpha =: (\gr^F A)_{+, \alpha}^{\otimes n}$.
\begin{lemma} The map 
\begin{equation*}
\varphi: (\gr^F A)_{+, \alpha}^{\otimes n} \rightarrow \frac{F_\alpha A_+^{\otimes n}}{F_{s(\alpha, -1)} A_+^{\otimes n}}
\end{equation*}
via 
\begin{equation*}
\varphi((a_1 + F_{s(\alpha_1,-1)} A) \otimes \cdots \otimes (a_n + F_{s(\alpha_n,-1)} A)) := 
a_1 \otimes \cdots \otimes a_n + F_{s(\alpha, -1)} A_+^{\otimes n}
\end{equation*}
is a chain isomorphism.
\end{lemma}

\begin{proof}
First, if $a_i - a'_i \in F_{s(\alpha_i, -1)} A$ for some $1 \leq i \leq n$ and $\alpha_1 \cdots \alpha_n = \alpha$, then 
\begin{equation*}
a_1 \otimes \cdots \otimes (a_i - a'_i) \otimes \cdots \otimes a_n \in F_{s(\alpha, -1)} A_+^{\otimes n}.
\end{equation*}
Hence, $\varphi$ is well-defined.

To show that $\varphi$ is a chain map, suppose $a_i \in F_{\alpha_i} A$ and $\alpha_1 \cdots \alpha_n = \alpha$.
We compute
\begin{align*}
\begin{split}
(d \circ \varphi)&\left(\left(a_1 + F_{s(\alpha_1, -1)} A\right) \otimes \cdots \left(a_n + F_{s(\alpha_n, -1)} A\right)\right)\\ 
&= d\left(a_1 \otimes \cdots \otimes a_n + F_{s(\alpha, -1)} A_+^{\otimes n}\right)
\end{split}
\\
&= \sum_{i = 1}^{n-1} (-1)^i a_1 \otimes \cdots \otimes a_i a_{i+1} \otimes \cdots \otimes a_n + F_{s(\alpha, -1)} A_+^{\otimes n}\\ 
\begin{split}
& = \varphi\Big(\sum_{i = 1}^{n-1} (-1)^i \left(a_1 + F_{s(\alpha_1, -1)} A\right) \otimes \cdots  \\
&  \otimes \left(a_i a_{i+1} + F_{s(\alpha_i \alpha_{i+1}, -1)} A\right) \otimes 
\cdots \otimes \left(a_n + F_{s(\alpha_n, -1)} A \right) \Big)
\end{split}
\\
& = (\varphi \circ d) \left(\left(a_1 + F_{s(\alpha_1, -1)} A\right) \otimes \cdots 
\otimes \left(a_n + F_{s(\alpha_n, -1)} A\right)\right).
\end{align*}

Now, to show that $\varphi$ is an isomorphism, note that the set 
\begin{equation*}
\B_1  := \Big\{\left(a_1 + F_{s(\alpha_1, -1)} A\right) \otimes \cdots \otimes \left(a_n + F_{s(\alpha_n, -1)} A\right) \Big|\;
 a_i \in \R_{\alpha_i}, \alpha_i \not = e, \alpha_1 \cdots \alpha_n = \alpha \Big\}
\end{equation*}
is a basis for $(\gr^F A)_{+, \alpha}^{\otimes n}$, while
\begin{equation*} \label{e:b2def}
\B_2  := \Big\{a_1 \otimes \cdots \otimes a_n + F_{s(\alpha, -1)} A_+^{\otimes n} \Big|\;
 a_i \in \R_{\alpha_i}, \alpha_i \not = e, \alpha_1 \cdots \alpha_n = \alpha \Big\}
\end{equation*}
is a basis for $F_\alpha A_+^{\otimes n}/F_{s(\alpha, -1)} A_+^{\otimes n}$.
Since $\varphi$ gives a bijection between these bases, $\varphi$ is an isomorphism.
\end{proof}

Now, because of condition (\ref{i:filtrationdimcondo}) on the filtration, we have a chain isomorphism
\begin{equation*}
\varphi^\vee : \left(\frac{F_\alpha A_+^{\otimes n}}{F_{s(\alpha, -1)} A_+^{\otimes n}}\right)^\vee \stackrel{\sim}{\rightarrow} \Cob_{\Gr}^{n, \alpha}(\gr^F A).
\end{equation*}
The restriction map 
\begin{equation*}
(A_+^{\otimes n})^\vee \rightarrow (F_\alpha A_+^{\otimes n})^\vee
\end{equation*}
induces an injective map
\begin{equation*}
\rho: \frac{F_\alpha \Cob^n(A)}{F_{s(\alpha, 1)} \Cob^n(A)} 
\hookrightarrow
\left(\frac{F_\alpha A_+^{\otimes n}}{F_{s(\alpha, -1)} A_+^{\otimes n}}\right)^\vee
.\end{equation*}
It is straightforward to check the following:
\begin{lemma} The map $\rho$ is a chain isomorphism.
\end{lemma}

We now know that 
\begin{equation*}
\varphi^\vee \circ \rho: \gr^F \Cob^\bullet(A) \rightarrow \Cob_{\Gr}^\bullet(\gr^F A)
\end{equation*}
is a chain isomorphism, graded by $\M$.

\begin{proof}[Proof of Proposition \ref{t:dgaisom}]
It suffices to show that $\varphi^\vee \circ \rho$ is a differential-graded algebra homomorphism. Let $f \in F_\alpha \Cob^n(A)$, $g \in F_\beta \Cob^m(A)$, $a_i \in F_{\alpha_i} A$, $b_i \in F_{\beta_i} A$, $\alpha_1 \cdots \alpha_n = \alpha$, and $\beta_1 \cdots \beta_m = \beta$.

Then, 
\begin{align*}
\begin{split}
(\varphi^\vee \circ \rho) & ((f + F_{s(\alpha, 1)} \Cob^n(A)) \smallsmile (g + F_{s(\beta, 1)}\Cob^m(A))\\
& \Big((a_1 + F_{s(\alpha_1, -1)} A) \otimes \cdots \otimes (a_n+ F_{s(\alpha_n, -1)} A)\\ 
&\otimes (b_1 + F_{s(\beta_1, -1)} A) \otimes \cdots \otimes (b_m + F_{s(\beta_m, -1)} A)\Big) \\
& = \rho((f + F_{s(\alpha, 1)} \Cob^n(A)) \smallsmile (g + F_{s(\beta, 1)} \Cob^m(A))) \\
& ((a_1 \otimes \cdots \otimes a_n + F_{s(\alpha, -1)} A_+^{\otimes n}) \otimes (b_1 \otimes \cdots \otimes b_m + F_{s(\beta, -1)} A_+^{\otimes m}))
\end{split}\\
\begin{split}
& = \rho(f \smallsmile g + F_{s(\alpha\beta, 1)}\Cob^{n+m}(A)) \\
& (a_1 \otimes \cdots \otimes a_n \otimes b_1 \otimes \cdots \otimes b_m + F_{s(\alpha\beta, -1)} A_+^{\otimes n + m})
\end{split}\\
& = f(a_1 \otimes \cdots \otimes a_n)g(b_1 \otimes \cdots \otimes b_m).
\end{align*}

Likewise,
\begin{align*}
\begin{split}
\Big((\varphi^\vee \circ \rho)&(f + F_{s(\alpha, 1)}\Cob^n(A)) \smallsmile (\varphi^\vee \circ \rho)(g + F_{s(\beta, 1)} \Cob^m(A))\Big)\\
& \Big((a_1 + F_{s(\alpha_1, -1)} A) \otimes \cdots \otimes (a_n + F_{s(\alpha_n, -1)} A)\\
& \otimes (b_1 + F_{s(\beta_1, -1)} A) \otimes \cdots \otimes (b_m + F_{s(\beta_m, -1)} A)\Big)\\
& = \rho(f + F_{s(\alpha, 1)} \Cob^n(A))(a_1 \otimes \cdots \otimes a_n + F_{s(\alpha, -1)} A_+^{\otimes n}) \\
& \cdot \rho(g + F_{s(\beta, 1)} \Cob^m(A))(b_1 \otimes \cdots \otimes b_m + F_{s(\beta, -1)} A_+^{\otimes m})
\end{split}\\
& = f(a_1 \otimes \cdots \otimes a_n)g(b_1 \otimes \cdots \otimes b_m),
\end{align*}
as desired.
\end{proof}

Recall that we give $\E(A)$ a filtration $F_\alpha \E(A)$ induced by the filtration $F_\alpha \Cob^\bullet(A)$.

\begin{definition}
Define a surjective map
 $\eta_\infty: F_\alpha \Cob^n(A) \cap \ker \partial \twoheadrightarrow (\gr^F \E(A))^{n, \alpha}$ to be the composition
\begin{equation*}
F_\alpha \Cob^n(A) \cap \ker \partial \twoheadrightarrow F_\alpha \E^n(A) \twoheadrightarrow \frac{F_\alpha \E^n(A)}{F_{s(\alpha, 1)} \E^n(A)}.
\end{equation*}
Define a map $\eta_1: F_\alpha \Cob^n(A) \cap \ker \partial \rightarrow \EGr^{n, \alpha}(\gr A)$ to be the composition
\begin{equation*}
\begin{split}
F_\alpha \Cob^n(A) \cap \ker \partial & \rightarrow 
\frac{F_\alpha \Cob^n(A) \cap \ker \partial + F_{s(\alpha, 1)}\Cob^n(A)}
{F_{s(\alpha, 1)} \Cob^n(A)} \\[1 em]
& \rightarrow \frac{F_\alpha \Cob^n(A)}{F_{s(\alpha, 1)} \Cob^n(A)} \cap \ker (\gr^F \partial)\\[1 em]
& \xrightarrow{\left. \varphi^\vee \circ \rho\right|} \Cob_{\Gr}^{n, \alpha}(\gr^F A) \cap \ker \partial \\[1 em]
& \twoheadrightarrow \EGr^{n, \alpha}(\gr^F A).
\end{split}
\end{equation*}
(Recall that $\varphi^\vee \circ \rho: \gr^F \Cob^\bullet(A) \rightarrow \Cob^\bullet_{\Gr}(\gr^F A)$ is a differential-graded algebra isomorphism by Proposition \ref{t:dgaisom}.)
\end{definition}
The maps $\eta_1$ and $\eta_\infty$ appear in the construction of a spectral sequence obtained from the filtration $F$ on $\Cob(A)$. See, for example, \cite[Theorem 2.6]{mcclearyspectral} and its proof. (We will not need this spectral sequence.)

\begin{lemma}\label{t:lambdawelldef} $\ker \eta_1 = \ker \eta_\infty$. 
\end{lemma}

\begin{proof}
Suppose $f \in \ker \eta_1$, meaning 
\begin{equation*}
(\varphi^\vee \circ \rho)(f + F_{s(\alpha, 1)} \Cob^n(A)) \in \Cob_{\Gr}^{n, \alpha}(\gr^F A) \cap \img \partial.
\end{equation*}
As $\varphi^\vee \circ \rho$ is a differential-graded algebra isomorphism, 
\begin{equation*}
f + F_{s(\alpha, 1)} \Cob^n(A) \in \frac{F_\alpha \Cob^n(A)}{F_{s(\alpha, 1)} \Cob^n(A)} \cap \img \partial;
\end{equation*} 
that is, there exists $g \in F_\alpha \Cob^{n-1}(A)$ such that
\begin{equation*}
\partial(g) + F_{s(\alpha, 1)} \Cob^n(A) = f + F_{s(\alpha, 1)} \Cob^n(A).
\end{equation*}
However, $f - \partial(g) + \img \partial \in F_{s(\alpha, 1)} \E^n(A)$. Thus, $\eta_\infty(f) = 0$.

Now, suppose $f \in \ker \eta_\infty$, meaning $f + \img \partial \in F_{s(\alpha, 1) }\E^n(A)$. So, $f + \partial(g) \in F_{s(\alpha,1)} \Cob^n(A)$ for some $g \in \Cob^n(A)$.
Since $f, f+ \partial(g) \in F_\alpha \Cob^n(A)$, $\partial(g) \in F_\alpha \Cob^n(A)$ as well, and
\begin{equation*}
f + F_{s(\alpha, 1)} \Cob^n(A) = \partial(g) + F_{s(\alpha, 1)} \Cob^n(A).
\end{equation*}
Thus,
\begin{equation*}
(\varphi^\vee \circ \rho)(f + F_{s(\alpha, 1)} \Cob^n(A)) \in \Cob_{\Gr}^{n, \alpha}(\gr ^F A) \cap \img \partial
\end{equation*}
and so $\eta_1(f) = 0$.
\end{proof}

\begin{definition} 
Since $\eta_\infty$ is surjective, Lemma \ref{t:lambdawelldef} tells us we may define a unique injective map $\Lambda^{n, \alpha}$ such that the diagram
\begin{equation*}
\xymatrix{ & F_\alpha \Cob^n(A) \cap \ker \partial \ar[ddl]_{\eta_\infty} \ar[ddr]^{\eta_1}& \\
\\
(\gr^F \E(A))^{n, \alpha} \ar[rr]^{\Lambda^{n, \alpha}} && \EGr^{n, \alpha}(\gr^F A)}
\end{equation*}
commutes. Set $\Lambda := \bigoplus_{n, \alpha} \Lambda^{n, \alpha}$.
\end{definition}

We may now prove Theorem \ref{t:algmonomorphism}, which we restate: 
\begin{theorem} 
The map
\begin{equation*}
\Lambda: \gr^F \E(A) \hookrightarrow \EGr(\gr^F A).
\end{equation*}
is an algebra monomorphism.
\end{theorem}

\begin{proof}
It remains only to prove $\Lambda$ is an algebra homomorphism.
Let $f \in (\gr^F \E(A))^{n, \alpha}$ and $g \in (\gr^F \E(A))^{m, \beta}$.
Choose preimages (under $\eta_1$) 
\begin{equation*}
\tilde{f} \in F_{\alpha} \Cob^n(A) \cap \ker \partial 
\text{ and } \tilde{g} \in F_{\beta} \Cob^m(A) \cap \ker \partial
\end{equation*}
 for $f$ and $g$, respectively. We have 
\begin{equation*}
\tilde{f} \otimes \tilde{g} \in F_{\alpha\beta}\Cob^{n + m}(A) \cap \ker \partial.\end{equation*}

Now, we compute
\begin{align*}
\eta_\infty(\tilde{f} \otimes \tilde{g}) & = ((\tilde{f} \otimes \tilde{g}) + \img \partial) + F_{s(\alpha\beta, 1)} \E^{n + m}(A) \\
& = ((\tilde{f} + \img \partial) \smallsmile (\tilde{g} + \img \partial)) + F_{s(\alpha\beta, 1)} \E^{n + m}(A)\\
& = ((\tilde{f} + \img \partial) + F_{s(\alpha, 1)} \E^n(A)) \smallsmile ((\tilde{g} + \img \partial) + F_{s(\beta, 1)} \E^m (A))\\
& = \eta_\infty(\tilde{f}) \smallsmile \eta_\infty(\tilde{g})\\
& = f \smallsmile g.\qedhere
\end{align*}
\end{proof}

Before proving Theorem \ref{t:kmpush}, we prove a general fact about filtered algebras:

\begin{lemma}\label{t:generalfactfiltered}
Let $R = \bigoplus_i R_i$ be a graded algebra with a decreasing filtration $F$ by a totally-ordered monoid $\M$. Put $F_\alpha R_i = F_\alpha R \cap R_i$ and assume $F_\alpha R = \bigoplus_i F_\alpha R_i$ for all $i$. Let $R'$ be the subalgebra of $R$ generated by $R_1, \dots, R_m$. Suppose, for each $i$, $F_\alpha R_i \subset R'$ for $\alpha$ sufficiently large. If  $(\gr^F R)_1, \dots (\gr^F R)_m$ generate $\gr^F R$, then $R_1, \dots, R_m$ generate $R$.
\end{lemma}

\begin{proof}
Suppose that $F_{s(\alpha, 1)} R_i \subset R'$. Let $a \in F_\alpha R_i \setminus F_{s(\alpha,1)} R_i$. As $\gr^F R$ is generated by $(\gr^F R)_1, \dots, (\gr^F R)_n$, there exists $a' \in R' \cap F_\alpha R_i$ such that
\begin{equation*}
a - a' \in F_{s(\alpha, 1)} R_i \subset R'.
\end{equation*} 
As $a \in R'$, we know $a \in R'$. Thus, $F_{\alpha} R \subset R'$. By (decreasing) induction on $\alpha$, $R = R'$.
\end{proof}

\begin{lemma} \label{t:findimgr}
If $\dim \gr^F \E^n(A) < \infty$ then $F_\alpha \E^n(A) = 0$ for some $\alpha$, and consequently, $\dim \gr^F \E^n(A) = \dim \E^n(A)$
\end{lemma}

\begin{proof}
Let $\{\xi + F_{s(\alpha_i, 1)} \E^n(A): 1 \leq i \leq m\}$ be a basis for $\gr^F \E^n(A)$, and choose $\alpha > \alpha_i$ for all $1 \leq i \leq m$. For $\beta \geq \alpha$, $F_\beta \E^n(A)/ F_{s(\beta, 1)} \E^n(A) = 0$, meaning $F_\beta \E^n(A) = F_\alpha \E^n(A)$.

Now, choose any $\xi \in F_\alpha \E^n(A)$. For $\beta \geq \alpha$, there exists $f_\beta \in F_\beta \Cob^n(A)$ and $f'_\beta \in \Cob^{n-1}(A)$ such that $f_\beta + \img d = \xi$ and $f_\beta = f_{s(\beta,1)} + d(f'_\beta)$.

Then, for $\beta \geq \alpha$,
\begin{align*}
f_\alpha &= f_{s(\alpha,1)} + d(f'_\alpha)\\
&= f_{s(\alpha,2)} + d(f'_{s(\alpha, 1)}) + d(f'_\alpha)\\
& \text{\hskip 1 cm} \vdots\\
& = f_\beta + \sum_{\alpha \leq \gamma < \beta} d(f'_\gamma).
\end{align*}
So, for $x \in F_\beta A_+^{\otimes n}$ and $\gamma > \beta$,
\begin{align*}
f_\alpha(x) &= f_\gamma(x) + \sum_{\alpha \leq \delta < \gamma} d(f'_\gamma)(x)\\
& = \sum_{\alpha \leq \delta < \gamma} (f'_\delta \circ \partial)(x).
\end{align*}
Thus, there exists $f':A_+^{\otimes n-1} \rightarrow \bbk$ such that $f_\alpha = f' \circ \partial$.  
Therefore, $\xi = 0$.
\end{proof}

We may now prove Theorem \ref{t:kmpush}, which we restate:

\begin{theorem}
If $\EGr^1(\gr^F A)$ and $\EGr^2(\gr^F A)$ are finite dimensional and generate $\EGr(\gr^F A)$, and $\Lambda^1$ and $\Lambda^2$ are surjective, then $A$ is $\ktwo$.
\end{theorem}

\begin{proof}
The map $\Lambda$ is an algebra isomorphism. Apply Lemma \ref{t:generalfactfiltered} when $m = 2$ and $R = \E(A)$.
\end{proof}

\begin{example}
Let \begin{equation*}
A = \frac{\bbk[x, y]}{\left<x^3 - p\right>},
\end{equation*}
where $p$ is a homogeneous quadratic polynomial. 
Define $\varepsilon: A \twoheadrightarrow \bbk$ via $\varepsilon(x) := 0$ and $\varepsilon(y) := 0$.

The standard $\mathbb{N}$-grading on $\bbk\left<x,y\right>$ induces a filtration $F$ on $A$  which satisfies the conditions in the introduction. 
Then,
\begin{equation*}
\gr^F A \simeq \frac{\bbk[x,y]}{\left<x^3\right>}.
\end{equation*}
Note that $\gr^F A$ is a complete intersection, and therefore is $\ktwo$ by \cite[Corollary 9.2]{cassidysheltonK2alg06}.

 One can easily compute $\dim E^1(\gr^F A) = \dim E^2(\gr^F A) = 2$. Furthermore, using $\Cob^\bullet(A)$, one can find the necessary linearly-independent cohomology classes to show $\dim \E^1(A) = \dim \E^2(A) = 2$, implying that $\Lambda^1$ and $\Lambda^2$ are surjective. Hence $A$ is $\ktwo$.
\end{example}


\section{Connected-graded algebras with monomial filtrations}
\label{s:monofilt}

By a {\bf connected-graded algebra}, we mean an algebra $A$ such that there is a graded algebra epimorphism
\begin{equation*}
\pi: \T(V) \twoheadrightarrow A
\end{equation*}
where $V = \spn\{x_1, \dots, x_n\}$ and $I := \ker \pi \subset \sum_{n \geq 2} V^{\otimes n}$ is finitely-generated and homogeneous. Under these circumstances, $\E^1(A)$ and $\E^2(A)$ are finite dimensional. The following lemma shows that the two definitions of $\ktwo$ from the introduction are compatible for connected-graded algebras.

\begin{lemma}\label{t:findimext}
For a connected-graded algebra $A$, $\E^m(A) = \EGr^m(A)$ if an only if $\dim \EGr^m(A) < \infty$. Consequently, if a connected-graded algebra $A$ is $\ktwo$ in the sense of Definition \ref{d:ktwoungraded}, then $A$ is $\ktwo$ in the sense of Definition \ref{d:ktwograded}.
\end{lemma}

\begin{proof}
Projective modules in the category $\Gr$-$A$ of graded $A$-modules are graded-free \cite[Proposition 2.1]{berger01}. So, there exists a projective resolution (in both the category of graded $A$-modules and of all $A$-modules) 
\begin{equation*}
\cdots \rightarrow A \otimes V^m \xrightarrow{\partial^m} \cdots \rightarrow A \otimes V^1 \rightarrow A \otimes V^0 \rightarrow A \rightarrow _A\bbk \rightarrow 0
\end{equation*}
such that each $V^i$ is a graded vector space and $\partial^i(A \otimes V^i) \subseteq A_+ \otimes V^{i-1}$. So, for \emph{any} $A$-module homomorphism $f: A \otimes V^{i-1} \rightarrow \bbk$, $\partial^i \circ f = 0$. Thus, all the differentials in both $\Hom(A \otimes V^\bullet, {}_A\bbk)$ and $\Hom_{\Gr}(A \otimes V^\bullet, {}_A\bbk)$ are zero. So, 
\begin{equation*}
E^m(A) = \Hom(A \otimes V^m,{}_A\bbk) \text{ while } \EGr^m(A) = \Hom_{\Gr}(A \otimes V^m, {}_A\bbk).\qedhere
\end{equation*}
\end{proof}

For a graded algebra $A$, we use notation established by \cite{cassidysheltonK2alg06}, setting
\begin{equation*}
A(n_1^{j_1}, n_2^{j_2}, \dots, n_t^{j_t}) := \bigoplus_i A(n_i)^{\oplus j_i}.
\end{equation*}

\begin{example} 
Consider the algebra 
\begin{equation*}
A = \frac{\bbk\left<w, x, y, z\right>}{\left<yz, zx-xz, zw\right>}.
\end{equation*}
introduced in \cite[Example 5.2]{cassidysheltonPBWdef}. A minimal projective resolution for $_A\bbk$ is
\begin{equation*}
0 \rightarrow A(-3, -4, -5, \dots) \rightarrow A(-2^2) \rightarrow A(-1^4) \rightarrow A \rightarrow \bbk \rightarrow 0.
\end{equation*}
Thus, the dimension of $\EGr^3(A)$ is countably infinite, while the dimension of $\E^3(A)$ is uncountable.
\end{example}

In light of Lemma \ref{t:findimext}, we will write $\E^1(A)$ for $\EGr^1(A)$ and $\E^2(A)$ for $\EGr^2(A)$.

The monomials $\M$ of $\T(V)$ (with respect to the basis $\{x_1, \dots, x_n\}$ for $V$) form a monoid which is totally-ordered by degree-lexicographical order. For $\alpha \in \M$, we set $F_\alpha A := \spn\left\{\pi(\beta): \beta \leq \alpha\right\}$. As $\M$ is itself $\mathbb{N}$-graded, the we may put an $\mathbb{N}$-grading on $\EGr(\gr^F A)$ by setting
\begin{equation*}
\EGr^{i,j}(\gr^F A) := \bigoplus_{|\alpha| = j}\EGr^{i, \alpha}(\gr^F A).
\end{equation*}
The algebra $\E(A)$ inherits the grading on $A$, and so does $\gr^F E(A)$. Indeed, it is clear that
\begin{equation*}
(\gr^F \E(A))^{i,j} = \bigoplus_{|\alpha| = j} (\gr^F \E^j(A))^{\alpha}
\end{equation*}
Furthermore, the  monomorphism \begin{equation*}
\Lambda: \gr^F \E(A) \hookrightarrow \EGr(\gr^F A)
\end{equation*} 
defined in Theorem \ref{t:algmonomorphism} is homogeneous with respect to this internal $\mathbb{N}$-grading. 

The goal of this section is to apply Theorem \ref{t:kmpush} to connected-graded algebras, using this monomial filtration. Note that $\Lambda^1$ is always surjective, so to apply Theorem \ref{t:kmpush}, we need only check:
\begin{enumerate}
\item $\gr^F A$ is $\ktwo$, and
\item $\Lambda^2: \gr^F \E^2(A) \hookrightarrow \E^2(\gr^F A)$ is surjective.
\end{enumerate}

Fortunately, the first condition is very easy to check since, as we will see, $\gr^F A$ is a monomial algebra.
\begin{definition}
\begin{enumerate}
\item We can write any element $x \in \T(V)$ uniquely as $c_\alpha \alpha + \sum_{\beta < \alpha} c_\beta \beta$ where $c_\alpha \not = 0$. Let $\tau(x) := c_\alpha \alpha$, which we call the {\bf leading monomial} of $x$. 
\item Define $\hat \pi: \T(V) \rightarrow \gr^F A$ via $\hat \pi(\alpha) = \pi(\alpha) + F_{s(\alpha, -1)} A$. 
\end{enumerate}
\end{definition}
We shall omit the proof of the following lemma.
\begin{lemma}[{\cite[Theorem 2.1]{math.RA/0609583}}] \label{t:jgendbytau}
$\ker(\hat \pi) = \left<\tau(x) : x \in I\right>$.
\end{lemma} 
From this, we see that $\gr^F A$ is a monomial algebra. Cassidy and Shelton
provide an algorithm that determines exactly when a monomial algebra is $\ktwo$
\cite[Theorem 5.3]{cassidysheltonK2alg06}. 

Now, we turn our attention to the second condition, the surjectivity of 
 $\Lambda^2: \gr^F \E^2(A) \hookrightarrow \E^2(\gr^F A)$. Let $I' = V \otimes I + I \otimes V$. 

\begin{definition}\label{d:essgrob}
\begin{enumerate}
\item An element $x \in I$ is an {\bf essential relation} for $A$ if $x$ is homogeneous and $x \not \in I'$. 
\item A generating set $\B^e$ for $I$ is an {\bf essential generating set} for $I$ if $\B^e$ comprises only essential relations and no subset of $\B^e$ generates $I$.
\item A generating set $\G$ for $I$ is a {\bf Gr\"obner basis} for $I$ if 
\begin{equation*}
\left<\tau(x): x \in I\right> = \left<\tau(x) : x \in \G\right>.
\end{equation*}
\item A Gr\"obner basis $\G$ for $I$ is an {\bf essential Gr\"obner basis} for $I$ if it is an essential generating set.
\end{enumerate}
\end{definition}

The definition of an essential relation first appeared in \cite{cassidysheltonK2alg06}. Gr\"obner bases are studied extensively in \cite{math.RA/0609583} and \cite{linoncomm}. 

Note that a generating set $\B^e$ for $I$ is essential if and only if $|\B^e| = \dim I/I' = \dim E^2(A)$. We will show later that the existence of an essential Gr\"obner basis is equivalent to the surjectivity of $\Lambda^2$. At the same time, it is desirable to know when an essential generating set is a Gr\"obner basis. 

\begin{example}\label{e:lmp}
Consider the ideal $I := \left<x^3, y^2\right>$ in $\bbk\left<x,y\right>$. Under the order $x < y$, the set $\B^e := \{y^2, x^3 - y^2x \}$ is an essential generating set for $I$. However $\B^e$ is not a Gr\"obner basis. On the other hand, the slightly modified set $\mathcal{G} := \{y^2, x^3\}$ is an essential Gr\"obner basis. The failure of $I$ to be a Gr\"obner basis was due to the needless redundancy of leading monomials.
\end{example}

The following lemma is easy.

\begin{lemma}\label{t:leadingmonolemma}
Let $\B^e$ be an essential generating set for $I$. Then the following are equivalent:
\begin{enumerate}
\item $\tau(\B^e)$ is an essential generating set for $\left<\tau(\B^e)\right>$.
\item For every $r, r' \in \B^e$ and $\alpha', \alpha'' \in \M$, $\tau(r) \not \in \bbk \alpha' \tau(r') \alpha''$. 
\item For every $r, r' \in \B^e$ and $\alpha', \alpha'' \in \M$, $\tau(r) \not \in \bbk \tau(\alpha' r' \alpha'')$.
\end{enumerate}
\end{lemma}

\begin{definition} If an essential generating set $\B^e$ meets the equivalent conditions of Lemma \ref{t:leadingmonolemma}, we say $\B^e$ has the {\bf leading monomial property}.
\end{definition}

In Example \ref{e:lmp}, the set $\B^e$ failed to be a Gr\"obner basis because it failed to have the leading monomial property.

\begin{lemma} \label{l:esgrobsatisfies} Essential Gr\"obner bases have the leading monomial property.\end{lemma}

\begin{proof}
Suppose that $\G$ is an essential Gr\"obner basis.  As we have an injective map $\Lambda^2: \gr^F \E^2(A) \hookrightarrow \EGr^2(\gr^F A)$, 
\begin{equation*}
|\G| = \dim \E^2(A) \leq \dim \EGr^2(\gr^F A).
\end{equation*}

On the other hand, if $\tau(r) \in \bbk\tau(\alpha' r' \alpha'')$ for some $\alpha', \alpha'' \in \M$ and $r, r' \in \B^e$, then 
\begin{equation*}
\left<\tau(\G)\right> = \left<\tau(\G) \setminus \{\tau(r)\}\right> 
\end{equation*}
and so $\dim \EGr^2(\gr^F A) < \dim \E^2(A)$, which is absurd.\end{proof}

\begin{theorem}\label{t:essentialspace} There exist homogeneous bases $\B$ for $I$ and $\B'$ for $I'$ such that $\B' \subset \B$, and the essential generating set $\B^e :=\B \setminus \B'$ has the leading monomial property.\end{theorem}

The proof of this theorem will follow after two technical lemmas.

\begin{lemma}\label{t:drilldownmonos} For $W \subset I$ and $\alpha \in \M$, define\label{ind:scriptA} 
\begin{equation*}
\mathcal{A}_m^\alpha(W) := \{r \in I_m : r \not \in \mathrm{span}\; W, \tau(r) \not \in \bbk \tau(s) \text{ for any } s \in W \text{ with } \tau(s) \geq \alpha\}.\end{equation*}
If $\mathcal{A}_m^{x_1^{m+1}}(W) \not = \emptyset$, then $\mathcal{A}_m^{x_1^m}(W) \not = \emptyset$; that is, there exists $r \in I_m$ such that $\tau(r) \not \in \bbk \alpha' \tau(s) \alpha''$ for any $\alpha', \alpha'' \in \mathcal{M}$ and $s \in W$.
\end{lemma}

\begin{proof}
Need only show that $\mathcal{A}_m^\alpha(W) \not = \emptyset$ implies that $\mathcal{A}_m^{s(\alpha,-1)}(W) \not = \emptyset$. Let $r \in \mathcal{A}_m^\alpha(W)$. Suppose $\tau(r) = \tau(s)$ for some $s \in W$. Then $r-s \in \mathcal{I}_m$ but $r-s \not \in \mathrm{span}\; W$. Also, $\tau(r-s) < \tau(s) < \alpha$, so $r-s \in \mathcal{A}_m^{s(\alpha, -1)}(W)$.\end{proof}

We will use the following lemma to build our basis degree-by-degree:
\begin{lemma}\label{t:degreeibuildess} Suppose $\B$ is a homogeneous basis for $\bigoplus_{i = 0}^{m-1} I_i$ and $\B'\subset \B$ is a basis for $\bigoplus_{i = 0}^{m-1} I'_i$. Then there exits $\B'' \subset I_m$ and $r_1, \dots, r_\ell \in I_m$ such that: 
\begin{enumerate}
\item $\B''$ is a basis for $I'_m$.
\item $r_i \not \in \bbk \alpha' \tau(r) \alpha''$ for any $i = 1, \dots, \ell$, 
$\alpha', \alpha'' \in \M$, and $r \in \B$.
\item $\B''\cup \{r_1, \dots, r_\ell\}$ is a basis of $I_m$.
\end{enumerate}
\end{lemma}

\begin{proof}
Set 
\begin{equation*}
\B^{(0)} = \left\{\alpha' r' \alpha'' \in I_m : \alpha', \alpha'' \in \mathcal{M}, r' \in \B\right\}.
\end{equation*}
Let $\B'' \subset \B^{(0)}$ such that $\B'_m$ is linearly independent. Since $\B^{(0)}$ spans $I'_m$, $\B''$ is a basis for $I'_m$.

Now, suppose we have constructed $\B^{(j)} = \B^{(j-1)} \cup \{r_j\}$ for $1 \leq j \leq i$ such that $(\B^{(i)} \setminus \B^{(0)}) \cup \B''$ is linearly independent and $\tau(r_j) \not \in \bbk\tau(s)$ for any $s \in \B^{(i-1)}$.

If $\B^{(i)}$ spans $I_m$, then 
$\B'' \cup \{r_1, \dots, r_i\}$ also spans $I_m$, and the claim is proved. Otherwise, $\mathcal{A}^{x_1^{m+1}}_m(\B^{(j)}) \not = \emptyset$, and so by Lemma \ref{t:drilldownmonos}, there exists $r_{i+1} \in I_m$ such that $\tau(r_{i+1}) \not \in \bbk \tau(s)$ for any $s \in \B^{(i)}$. Set $\B^{(i+1)} = \B^{(i)} \cup \{r_{i+1}\}$.
\end{proof}

\begin{proof}[Proof of Theorem \ref{t:essentialspace}] Set $\B_m = \B'_m = \B^e_m = \emptyset$ for $m \leq 1$. Apply Lemma \ref{t:degreeibuildess} and induction on $m$.
\end{proof}

We are now ready to prove Theorem \ref{t:genPBWa}, which we restate:
\vskip 5 mm 
\begin{theorem}\addtolength{\parskip}{-3mm} \label{t:genPBW}
The following are equivalent:
\begin{enumerate}
\item Every essential generating set for $I$ with the leading monomial property is a Gr\"obner basis. \label{t:genPBW:everygenset}
\item There is an essential Gr\"obner basis for $I$. \label{t:genPBW:existsgrob}
\item $\dim E^2(A) = \dim E^2(\gr^F A)$. \label{t:genPBW:dimsequal}
\item The injective map $\Lambda^2: \gr^F \E^2(A) \hookrightarrow \E^2(\gr^F A)$ defined in Theorem \ref{t:algmonomorphism} is surjective. \label{t:genPBW:lambdasurj}
\end{enumerate}
Therefore, if $I$ has an essential Gr\"obner basis and $\gr^F A$ is $\ktwo$, then $A$ is $\ktwo$ as well.
\end{theorem}

\begin{proof}
We set $J = \ker(\hat \pi: \T(V) \twoheadrightarrow \gr^F A)$ and $J' = J \otimes V + V \otimes J$.

In light of Theorem \ref{t:essentialspace}, is clear that condition (\ref{t:genPBW:everygenset}) implies condition (\ref{t:genPBW:existsgrob}). 

Suppose $\G$ is an essential Gr\"obner basis for $I$. Then $|\G| = \dim I/I'$. Also, since $\G$ has the leading monomial property, $|\G| = |\tau(\G)| = \dim J/J'$. So, condition (\ref{t:genPBW:existsgrob}) implies condition (\ref{t:genPBW:dimsequal}).

Clearly, condition (\ref{t:genPBW:dimsequal}) and condition (\ref{t:genPBW:lambdasurj}) are equivalent.

Finally, assume (\ref{t:genPBW:lambdasurj}). Suppose $\B^e$ is an essential generating set for $I$ with the leading monomial property. Let $\B^e_J$ be an essential generating set of $J$ such that 
\begin{equation*}
\{\tau(x): x \in \B^e \} \subset \B^e_J.
\end{equation*} 
Then, $|\B^e| = \dim I/I' = \dim J/J' = |\B^e_J|$. So, $\B^e$ is a Gr\"obner basis. Thus, condition (\ref{t:genPBW:lambdasurj}) implies condition (\ref{t:genPBW:everygenset}).
\end{proof}

\begin{example} \label{ex:assocgrows} 
Consider 
\begin{equation*}
A: = \frac{\bbk\left<x, y\right>}{\left<xy - x^2, yx, y^3\right>}
\end{equation*} 
with a monomial order induced by $x < y$. We know from \cite[Example 4.5]{cassidysheltonK2alg06} that $A$ is not a $\ktwo$ algebra. The Hilbert series of $A$ is $H_A(t) = 1 + 2t + 2t^2$. Since $\pi(x^3) = 0$, we see that $\hat\pi(x^3) = 0$, and $\gr^F A \simeq \bbk\left<x, y\right> / \left<xy, yx, x^3, y^3\right>$. We may apply \cite[Theorem 5.3]{cassidysheltonK2alg06} to see that $\gr^F A$ is $\ktwo$. The essential generating set $\{xy - x^2, yx, y^3\}$ is not a Gr\"obner basis for $\ker \pi$. The behavior is similar under $y < x$ (although $\gr^F A$ is a different $\ktwo$ algebra).
\end{example}

\begin{example} \label{ex:conversefalse}
Consider \begin{equation*}
A := \frac{\bbk\left<x, y, z\right>}{\left<x^2y - x^3, yz^2 - yx^2, x^3z- x^4\right>}
\end{equation*}
with the monomial order induced by $x < y < z$.
 We may use the diamond lemma \cite[Theorem 1.2]{bergman78} to show that 
\begin{equation*}
\gr^F A \simeq B:= \frac{\bbk\left<x, y, z\right>}{\left<x^2y, yz^2, x^3z\right>}.
\end{equation*} 
Thus,  $\{x^2y - x^3, yz^2 - yx^2, x^3z-x^4\}$ is an essential Gr\"obner basis for $\ker \pi$. However, application of \cite[Theorem 5.3]{cassidysheltonK2alg06} shows that $B$ is not $\ktwo$. By inspection, 
\begin{equation*}
0 \rightarrow B(-5) \xrightarrow{\left(\begin{matrix}0 & x^2 & 0\end{matrix}\right)}{\longrightarrow} B(-3^2, -4) \xrightarrow{\left(\begin{matrix} 0 & x^2 & 0 \cr
0 & 0 & yz \cr
0 & 0 & x^3
\end{matrix}\right)} B(-1^3) \xrightarrow{\left(\begin{matrix} x \cr y \cr z\end{matrix}\right)} B \rightarrow \bbk \rightarrow 0
\end{equation*}
is a minimal projective resolution for $_B\bbk$. By Theorem \ref{t:algmonomorphism}, $\dim \E^{i,j}(A) \leq \dim \E^{i,j}(B)$. So, the chain complex of projective $A$-modules
\begin{equation*}
0 \rightarrow  A(-5) \xrightarrow{\left(\begin{matrix}0 & x^2 & -x\end{matrix}\right)}
A(-3^2, -4) \xrightarrow{\left(\begin{matrix}
x^2 & -x^2 & 0 \cr
y^2 & 0 & -yx \cr
x^3 & 0 & -x^3
\end{matrix}\right)}
A(-1^3) \xrightarrow{\left(\begin{matrix} x \cr y \cr z \end{matrix}\right)} A \rightarrow \bbk \rightarrow 0
\end{equation*}
is a minimal projective resolution for $_A\bbk$. Applying \cite[Theorem 4.4]{cassidysheltonK2alg06}, we see that $A$ is $\ktwo$. Hence, the converse to Theorem \ref{t:genPBW} is false.
\end{example}

\begin{example}
Let 
\begin{equation*}
A:= \frac{\bbk\left<x,y\right>}{\left<yx -xy, y^3 + x^2y\right>}.
\end{equation*} 
Then under the order $x < y$, the essential generating set $\{yx-xy, y^3 + x^2y\}$ is a Gr\"obner basis for $\ker \pi$, and 
\begin{equation*}
\gr^F A = \frac{\bbk\left<x,y\right>}{\left<yx, y^3\right>}.
\end{equation*} 

We may use \cite[Theorem 5.3]{cassidysheltonK2alg06} to show that $\gr^F A$ is $\ktwo$. Thus, by Theorem \ref{t:genPBW}, $A$ is $\ktwo$. (This can also be verified directly using \cite[Corollary 9.2]{cassidysheltonK2alg06}.)
\end{example}

Theorem \ref{t:genPBW} is a generalization of the classical theory of Poincar\'e--Birkhoff--Witt algebras, which we can also prove:
\begin{theorem}[{\cite[Theorem IV.3.1]{polishquadalg}}]
If $A$ is a quadratic algebra, and $\gr^F A$ is also quadratic, then $A$ is Koszul.\end{theorem}
\begin{proof}
Quadratic monomial algebras are Koszul \cite[Corollary II.4.3]{polishquadalg}. The theorem follows directly from Theorem \ref{t:kmpush}.
\end{proof}


\section{Anticommutative analogues to face rings}
\label{s:anticommface}

In this section, use the results from Section \ref{s:monofilt} to show some anticommutative analogues to face rings are $\ktwo$. In particular, we prove:
\begin{theorem}\label{t:grassmannmodall}
 The algebra
\begin{equation*}
\frac{\bigwedge_\bbk (x_1, \dots, x_n)}{(x_1 \cdots x_n)}
\end{equation*}
is $\ktwo$.
\end{theorem}
Suppose $X := \{x_1, \dots, x_n\}$ is a finite set and $\Delta$ is a simplicial complex on $X$---that is, $\Delta \subset 2^X$ such that $\{x_i\} \in \Delta$ for $1 \leq i \leq n$ and if $Y \in \Delta$, then $2^Y \subset \Delta$. We define an algebra
\begin{equation*}
A[\Delta] := \bigwedge_\bbk(x_1, \dots, x_n)/\left<x_{i_1} \cdots x_{i_r}| i_1 < i_2 < \cdots < i_r, \{x_{i_1}, \dots x_{i_r}\} \not \in \Delta\right>,
\end{equation*}
the anticommutative analogue of the face ring of $\Delta$. (Face rings are studied in detail in \cite{stanleycombin}.)

\begin{definition}
If $Y \subset X$, $Y \not \in \Delta$, but $2^Y \setminus \{Y\} \subset \Delta$, then we say $Y$ is a {\bf minimally missing face} of $\Delta$. 
\end{definition}

\begin{theorem}\label{t:essentialgrobcomp} Suppose $\Delta$ is a simplicial complex on $X := \{x_1, \dots, x_n\}$. Under the order $x_1 < \dots < x_n$, $\ker \pi_{A[\Delta]}$ has an essential Gr\"obner basis if and only if every minimally missing face $Y:=\{x_{i_1}, \dots, x_{i_m}\} \subset X$ (where $i_1 < i_2 < \dots < i_m$) satisfies the following property: 
\begin{equation} \label{e:minmissingprop}\text{If $u \not \in Y$ and $i_1 < u <i_m$, then $(Y \setminus \{x_{i_1}\}) \cup \{x_u\} \not \in \Delta$ or $(Y \setminus \{x_{i_m}\}) \cup \{x_u\} \not \in \Delta$.}
\end{equation}
\end{theorem}

\begin{proof}
An essential generating set with the leading monomial property for
 \begin{equation*}
I := \ker(\pi: \bbk\left<x_1, \dots, x_n\right> \rightarrow A[\Delta])
\end{equation*}
 is
 \begin{equation*}
\begin{split}
\B^e = & \{x_jx_i + x_ix_j | i < j\} \cup \{x_i^2|i = 1 \dots n\} \\
& \cup \{x_{i_1} \cdots x_{i_m}  | i_1 < \dots < i_m, \{x_{i_1}, \dots, x_{i_m}\} \text{ is a minimally missing face}\}.
\end{split}
\end{equation*}

If $Y$ is a minimally missing face which fails (\ref{e:minmissingprop}) for some $u \not \in Y$, then 
\begin{equation*}
x_{i_1} \cdots x_{i_t} x_u x_{i_{t+1}} \cdots x_m
\end{equation*} 
is an essential relation of $\gr^F A$ for some $t$, meaning that 
$\B^e$ is not a Gr\"obner basis.

On the other hand, suppose $\B^e$ is not a Gr\"obner basis. 
Then $\gr^F A$ has some new essential relation $r$ such that $r \not = \tau(x)$ for $x \in \B^e$. Pick such $r$ minimally. Then 
\begin{equation*}
r = x_{i_1} \cdots x_{i_m}x_u \bmod \left<x_i x_j + x_j x_i\right>
\end{equation*} 
for some minimally missing face $Y = \{x_{i_1}, \dots, x_{i_m}\}$. So $Y$ fails (\ref{e:minmissingprop}).
\end{proof}

\begin{proof}[Proof of Theorem \ref{t:grassmannmodall}]
Let $X := \{x_1, \dots, x_n\}$ and $\Delta = 2^X \setminus \{X\}$.
Then by Theorem \ref{t:essentialgrobcomp}, $\ker(\pi: \bbk\left<x_1, \dots, x_n\right> \rightarrow A[\Delta])$ has an essential Gr\"obner basis.

So, applying \cite[Theorem 5.3]{cassidysheltonK2alg06} to
\begin{equation*}
\gr^F A = \bbk\left<x_1, \dots, x_n\right> / \left<x_1\cdots x_n, x_jx_i: 1 \leq i \leq j \leq n\right>,
\end{equation*}
we see that $\gr^F A$ is $\ktwo$, and hence $A$ is $\ktwo$.
\end{proof}

Not every simplicial complex $\Delta$ on a set $X$ has an ordering of $X$ which yields an essential Gr\"obner basis for $\ker \pi_{A[\Delta]}$.
\begin{example} Set $X := \{t, u, w, x, y, z\}$ and 
\begin{equation*}\begin{split}
\Delta & := \left(2^{\{u, x, y, z\}} \cup 2^{\{t, u, x, z\}} \cup 2^{\{u, w, x, z\}}\right) \setminus \{\{u,x,y,z\},\\ & \{t,u,x,z\}, \{u, w, x, z\}, \{x,y,z\}, \{t, u, z\}, \{u, w, x\}\}.
\end{split}
\end{equation*}

 Suppose we have an order $<$ of $X$ under which $\ker \pi_{A[\Delta]}$ has an essential Gr\"obner basis.

Note that $\{x, y, z\}$ is a minimally missing face, but $\{u, x, y\}, \{u, y, z\}, \{u, x, z\} \in \Delta$. So either $u < x,y,z$ or $u > x, y, z$. Without loss of generality, $u < x, y, z$.

Also, $\{t, u, z\}$ is a minimally missing face, but $\{u, x, z\}, \{t, x, z\}, \{t, u, x\} \in \Delta$. So as $u <x$, $x > t, u, z$.

Finally, $\{u, w, x\}$ is a minimally missing face, but $\{u, x, z\}, \{u, w, z\}, \{w, x, z\} \in \Delta$. However, as $x > z$, we cannot have $z > x, u, w$. However, as $u < z$, we cannot have $z < x, u, w$ either.
\end{example}

\bibliographystyle{amsplain}
\bibliography{genkoszul}

\providecommand{\bysame}{\leavevmode\hbox to3em{\hrulefill}\thinspace}
\providecommand{\MR}{\relax\ifhmode\unskip\space\fi MR }
\providecommand{\MRhref}[2]{%
  \href{http://www.ams.org/mathscinet-getitem?mr=#1}{#2}
}
\providecommand{\href}[2]{#2}
\begin{thebibliography}{10}

\bibitem{berger01}
Roland Berger, \emph{Koszulity for nonquadratic algebras}, Journal of Algebra
  \textbf{239} (2001), no.~2, 705--734. \MR{1832813 (2002d:16034)}

\bibitem{bergman78}
George~M. Bergman, \emph{The diamond lemma for ring theory}, Advances in
  Mathematics \textbf{29} (1978), no.~2, 178--218. \MR{506890 (81b:16001)}

\bibitem{cassidysheltonK2alg06}
Thomas Cassidy and Brad Shelton, \emph{Generalizing the notion of a {K}oszul
  algebra}, Mathematische {Z}eitschrift, to appear.

\bibitem{cassidysheltonPBWdef}
\bysame, \emph{{PBW}-deformation theory and regular central extensions},
  Journal f\"ur die reine und angewandte {M}athematik, to appear.

\bibitem{green04}
E.L. Green, E.~N. Marcos, R.~{Mart\'\i nez-Villa}, and Pu~Zhang,
  \emph{${D}$-{K}oszul algebras}, Journal of Pure and Applied Algebra
  \textbf{193} (2004), 141--162. \MR{2076383 (2005f:16044)}

\bibitem{linoncomm}
Huishi Li, \emph{Noncommutative {G}r\"obner bases and filtered-graded
  transfer}, Lecture notes in mathematics, vol. 1795, Springer, Berlin, 2002.
  \MR{1947291 (2003i:16065)}

\bibitem{math.RA/0609583}
\bysame, \emph{Algebras and their associated monomial algebras}, pre-print,
  {\tt arXiv:math/0609583v1 [math.RA]}, 2006.

\bibitem{mcclearyspectral}
John {McCleary}, \emph{A user's guide to spectral sequences}, second ed.,
  Cambridge studies in advanced mathematics, vol.~58, Cambridge University
  Press, Cambridge, 2001. \MR{1793722 (2002c:55027)}

\bibitem{polishquadalg}
Alexander Polishchuk and Leonid Positselski, \emph{Quadratic algebras},
  University lecture series, vol.~37, American Mathematical Society,
  Providence, RI, 2005. \MR{2177131 (2006f:16043)}

\bibitem{stanleycombin}
Richard~P. Stanley, \emph{Combinatorics and commutative algebra}, second ed.,
  Progress in mathematics, vol.~41, Birkh\"auser Boston, Boston, 1996.
  \MR{1453579 (98h:05001)}

\end{thebibliography}

\end{document}